\documentclass{article}

\usepackage{amsmath}
\usepackage{url}

\usepackage{amssymb}
\usepackage{amsthm}
\usepackage{graphicx}

\newcommand{\R}{\mathbb{R}}
\newcommand{\C}{\mathbb{C}}
\newcommand{\N}{\mathbb{N}}
\newcommand{\Q}{\mathbb{Q}}

\newcommand{\axright}{(\stackrel{\longrightarrow}{\text{S6}})}
\newcommand{\axleft}{(\stackrel{\longleftarrow}{\text{S6}})}
\newcommand{\sumrinline}{\mathop{\to\!\!\!\!\!\!\sum}}
\newcommand{\sumr}{\mathop{\to\!\!\!\!\!\!\!\!\sum}\limits}
\newcommand{\suml}{\mathop{\leftarrow\!\!\!\!\!\!\!\!\sum}\limits}

\newcommand{\prodr}{\mathop{\to\!\!\!\!\!\!\!\prod}\limits}

\def\lineclear
    {\rule{0pt}{0pt}\nopagebreak\par\nopagebreak\noindent}

\newtheorem{theorem}{Theorem}[section]
\newtheorem{proposition}[theorem]{Proposition}
\theoremstyle{definition}
\newtheorem{definition}[theorem]{Definition}

\begin{document}

\title{How to Add a Noninteger Number of Terms: From Axioms to New Identities}
\author{Markus M\"uller and Dierk Schleicher}
\date{}
\maketitle

\begin{abstract}
Starting from a small number of well-motivated axioms, we derive a unique definition of sums with a noninteger number of addends. These ``fractional sums'' have properties that generalize well-known classical sum identities in a natural way. We illustrate how fractional sums can be used to derive infinite sum and special functions identities; the corresponding proofs turn out to be particularly simple and intuitive.
\end{abstract}

\bigskip

\centerline{
``God made the integers; all else is the work of man.'' 
}
\hfill Leopold Kronecker

\section{Introduction.}
Mathematics is the art of abstraction and generalization.
Historically, ``numbers'' were first natural numbers; then
rational, negative, real, and complex
numbers were introduced (in some order). Similarly, the concept of taking derivatives has been generalized from first, second, and higher order derivatives to  ``fractional calculus'' of noninteger orders (see for instance \cite{FracCalculus}), and there is also some work on fractional iteration. 

However, when we add some number of terms, this
number (of terms) is still generally considered a natural number:
we can add two, seven, or possibly zero numbers, but what is the
sum of the first $-7$ natural numbers, or the first $\pi$ terms of the
harmonic series?

In this note, we show that there is a very natural way of
extending summations to the case when the ``number of terms'' is
real or even complex. One would think that this method should have been discovered at least two hundred years ago --- and that is what we initially suspected as well. To our surprise, this method does not seem to have been investigated in the literature, or to be known by the experts, apart from sporadic remarks even in Euler's work \cite{Euler} (see equation (\ref{Eq:Euler}) below). Of course, one of the standard methods to introduce the $\Gamma$ function is an example of a summation with a complex number of terms; we discuss this in Section~\ref{SecFromAxToDef}, equation (\ref{EqGamma}).\footnote{\textbf{Note by the second author.}
Many of the original ideas in this text are due to Markus M\"uller, who invented them while he was a high school student in the remote German province town of Morsbrunn, Bavaria. He was lacking the skills to carry out a formal mathematical proof, but he kept producing the most obscure mathematical identities on classical sums and fractional sums. I met him at the science contest ``Jugend forscht''  for high school students, and from then on helped him to turn his ideas into actual mathematical theorems and proofs, and to find out which of his formulas and identities were correct (most of them were). 
The main results have been published in \emph{The Ramanujan Journal} \cite{FracSums}.
}

Since this note is meant to be an introduction to
an unusual way of adding, we skip some of the proofs and
refer the reader instead to the more formal note \cite{FracSums}.
Some of our results were initially announced in \cite{FracSumsProc}.

\subsection{The Axioms.}
\label{SecTheAxioms}

We start by giving natural conditions for summations with
an arbitrary complex number of terms; here $x$, $y$, $z$, and $s$ are complex numbers and $f$ and $g$ are complex-valued functions defined on $\C$ or subsets thereof, subject to some conditions that we specify later:
\begin{description}
\item[(S1) Continued Summation]
\begin{equation}
\sum_{\nu=x}^y f(\nu)
+\sum_{\nu=y+1}^z f(\nu)
=\sum_{\nu=x}^z f(\nu).
\notag
\label{EqS1}
\end{equation}

\item[(S2) Translation Invariance]
\begin{equation}
\sum_{\nu=x+s}^{y+s} f(\nu)
=\sum_{\nu=x}^y f(\nu+s).
\notag
\label{EqS2}
\end{equation}

\item[(S3) Linearity] for arbitrary constants $\lambda,\mu\in\C$,
\begin{equation}
   \sum_{\nu=x}^y \left(\lambda f(\nu)+\mu g(\nu)\right) =
   \lambda\sum_{\nu=x}^y f(\nu)+\mu\sum_{\nu=x}^y g(\nu).
\notag
   \label{EqS3}
\end{equation}

\item[(S4) Consistency with Classical Definition]
\begin{equation}
   \sum_{\nu=1}^1 f(\nu)=f(1).
\notag
\label{EqS4}
\end{equation}

\item[(S5) Monomials]
for every $d\in\N$, the mapping
\begin{equation}
   z\mapsto\sum_{\nu=1}^z \nu^d
\notag
   \label{EqS6}
\end{equation}
is holomorphic in $\C$.

\item[\boldmath$\axright$ \unboldmath  Right Shift Continuity]
if $\lim_{n\to+\infty}f(z+n)=0$ pointwise for every $z\in\C$, then
\begin{equation}
   \lim_{n\to+\infty} \sum_{\nu=x}^y f(\nu+n)=0;
   \label{EqS5}
\end{equation}
more generally, if there is a sequence of polynomials $(p_n)_{n\in\N}$ of fixed degree 
such that, as $n\to +\infty$, $|f(z+n)-p_n(z+n)|\longrightarrow 0$ for all $z\in\C$, we require that
\begin{equation}
   \left\vert \sum_{\nu=x}^y f(\nu+n)-\sum_{\nu=x}^y p_n(\nu+n)\right\vert\longrightarrow 0.
\label{EqS5b}
\end{equation}

\end{description}

The first four axioms (S1)--(S4) are so obvious that it is hard to imagine any summation theory that violates these. They easily imply $\sum_{\nu=1}^n f(\nu)=f(1)+f(2)+\dots+f(n)$ for every $n\in\N$, so we are being consistent with the classical definition of summation.

Axiom (S5) is motivated by the well-known formulas 
\[
\sum_{\nu=1}^n \nu=\frac{n(n+1)}2,\quad
\sum_{\nu=1}^n \nu^2=\frac{n(n+1)(2n+1)}{6},\quad
\sum_{\nu=1}^n \nu^3 = \left(\frac{n(n+1)}{2}\right)^2
\]
and similarly for higher powers; we shall show below that our axioms imply that all those formulas remain valid for arbitrary $n\in\C$.

Finally, axiom $\axright$ is a natural condition also. The first case, in (\ref{EqS5}), expresses the view that if $f$ tends to zero, then the summation ``on the bounded domain'' $[x,y]$ should do the same. In (\ref{EqS5b}), the same holds, except an approximating polynomial is added; compare the discussion after Proposition~\ref{Prop1}. 

It will turn out that for a large class of functions $f$, there is
a unique way to define a sum $\sum_1^z f$ with $z\in\C$ that respects all
these axioms. In the next section, we will derive this definition and
denote such sums by $\sumrinline_1^z f$. We call them ``fractional sums.''

\subsection{From the Axioms to a Unique Definition.}
\label{SecFromAxToDef}

To see how these conditions determine a summation method
uniquely, we start by summing up polynomials. The simplest such
case is the sum $\sum_{\nu=1}^{\frac 1 2} c$ with $c\in\C$ constant.
If axiom (S1) is respected, then
\[
   \sum_{\nu=1}^{1/2} c +\sum_{\nu=3/2}^1 c =\sum_{\nu=1}^1 c.
\]
Applying axioms (S2) on the left and (S4) on the right-hand side, one gets
\[
   \sum_{\nu=1}^{1/2} c +\sum_{\nu=1}^{1/2} c =c.
\]
It follows that $\sum_{\nu=1}^{1/2} c =c/2$. This simple
calculation can be extended to cover every sum of polynomials with a rational number of terms.

\begin{proposition}
\label{Prop1}
For any polynomial $p:\C\to\C$, let $P:\C\to\C$ be the unique polynomial with $P(0)=0$ and $P(z)-P(z-1)=p(z)$ for all $z\in\C$.
Then:
\begin{itemize}
\item The possible definition
\begin{equation}
   \sum_{\nu=x}^y p(\nu):=P(y)-P(x-1)
   \label{eqSumPoly}
\end{equation}
satisfies all axioms (S1) to $\axright$ for the case that $f$ is a polynomial.
\item Conversely, every summation theory that satisfies axioms (S1), (S2), (S3), and (S4) also satisfies~(\ref{eqSumPoly})
for every polynomial $p$ and all $x,y\in\C$ with rational difference $y-x\in\Q$.
\item Every summation theory that satisfies (S1), (S2), (S3), (S4), and (S5) also satisfies~(\ref{eqSumPoly}) for every polynomial $p$ and all $x,y\in\C$.
\end{itemize}
\end{proposition}

\begin{proof}
To prove the first statement, suppose we use (\ref{eqSumPoly}) as a definition. It is trivial to check that this definition satisfies
(S1), (S3), (S4), and (S5). To see that it also satisfies (S2), consider a polynomial $p$ and the unique
corresponding polynomial $P$ with $P(x)-P(x-1)=p(x)$ and $P(0)=0$. 
Define $\tilde p(x):=p(x+s)$
and $\tilde P(x):=P(x+s)-P(s)$.
Then $\tilde P(0)=0$, and $\tilde P(x)-\tilde P(x-1)=P(x+s)-P(x+s-1)=p(x+s)=\tilde p(x)$. Hence
\[
   \sum_{\nu=x}^y p(\nu+s)=\sum_{\nu=x}^y\tilde p(\nu)=\tilde P(y)-\tilde P(x-1)=P(y+s)-P(x+s-1)=\sum_{\nu=x+s}^{y+s} p(\nu).
\]
To see that (\ref{eqSumPoly}) also satisfies $\axright$, let $V_\sigma$ be the linear space of complex
polynomials of degree less than or equal to $\sigma\in\N$. The definition
$\|p\|:=\sum_{i=0}^\sigma |p(i)|$ for $p\in V_\sigma$ introduces a norm on $V_\sigma$.
If we define a linear operator $\sum_x^y:V_\sigma\to\C$ via $\sum_x^y p := \sum_{\nu=x}^y p(\nu)$, then
this operator is bounded since $\dim V_\sigma=\sigma+1<\infty$. Thus, if $(q_n)_{n\in\N}\subset V_\sigma$ is
a sequence of polynomials with $\lim_{n\to\infty} \|q_n\|=0$, we have
$\lim_{n\to\infty} \sum_{\nu=x}^y q_n(\nu)=0$. Axiom $\axright$ then follows from considering
the sequence of polynomials $(q_n)_{n\in\N}$ with $q_n(x):=p(x+n)-p_n(x+n)$ and noting that
pointwise convergence to zero implies convergence to zero in the norm $\|\cdot\|$ of $q_n$, and thus of $\sum_x^y q_n$.

To prove the second statement,
we extend the idea that we used above to show that $\sum_{\nu=1}^{1/2} c=c/2$. Using (S1), we write for an integer $r\ge 1$
\[
\sum_{\nu=1}^{r}\nu^d
=
\sum_{\nu=1}^{r/s}\nu^d + \sum_{\nu=r/s+1}^{2r/s} \nu^d + \cdots + \sum_{\nu=(s-1)r/s+1}^{r}\nu^d,  
\]
where the left-hand side has a classical interpretation, using (S1), (S2), and (S4). Rewriting the right-hand side according to (S2)
and using (S3), we get
\[
\sum_{\nu=1}^{r}\nu^d
=
   \sum_{k=0}^{s-1} \sum_{\nu=1}^{r/s} \left(\nu+\frac{kr}s\right)^d=s\cdot \sum_{\nu=1}^{r/s} \nu^d + \sum_{k=0}^{s-1}
   \sum_{\nu=1}^{r/s} q_{d-1,k}(\nu),
\]
where the $q_{d-1,k}(\nu)=(\nu+kr/s)^d-\nu^d$ are polynomials of degree $d-1$ (and all $q_{-1,k}\equiv 0$). Now we argue by induction. If $d=0$,
the previous equation clearly determines $\sum_{\nu=1}^{r/s}1$ and by linearity also the corresponding sum over
arbitrary constants $c\in\C$. Once the value of the sum of any polynomial of degree $d-1$ is determined, the
equality also determines the value of $\sum_{\nu=1}^{r/s} \nu^{d}$, and by linearity, the sum of
every polynomial $p$ of degree $d$. Using (S2) again, we see that the axioms (S1) to (S4) uniquely determine
the sum $\sum_{\nu=x}^y p(\nu)=\sum_{\nu=1}^{y-x+1} p(\nu+x-1)$
if $y-x\in\Q$. As we have seen, equation~(\ref{eqSumPoly}) is a possible definition satisfying those axioms; hence it is the only
possible definition for $y-x\in\Q$.

Finally, it is clear how the restriction $y-x\in\Q$ can be lifted by additionally assuming (S5) (equation (\ref{eqSumPoly}) is already satisfied for $x-y\in\R$ by requiring just continuity in (S5); holomorphy is required for $x-y\in\C$).
\end{proof}

Consider now an arbitrary function $f\colon\C\to\C$. If we are
interested in $\sum_{\nu=x}^y f(\nu)$ for complex $x,y\in\C$, we can write
\[
   \sum_{\nu=x}^y f(\nu)+\sum_{\nu=y+1}^{y+n} f(\nu)=
   \sum_{\nu=x}^{x+n-1} f(\nu)+\sum_{\nu=x+n}^{y+n} f(\nu),
\]
where $n\in\N$ is an arbitrary natural number. Hence,
\begin{eqnarray}
   \sum_{\nu=x}^y f(\nu)
   &=&
   \sum_{\nu=x}^{x+n-1} f(\nu)
   -\sum_{\nu=y+1}^{y+n} f(\nu)
   +\sum_{\nu=x+n}^{y+n} f(\nu)
   \nonumber
   \\
   &=&
   \sum_{\nu=1}^{n} \left(f(\nu+x-1)-f(\nu+y)\strut\right)
   +\sum_{\nu=x}^{y} f(\nu+n).
   \label{EqHeuristics}
\end{eqnarray}
What have we achieved by this elementary rearrangement? In the last line, the first
sum on the right-hand side involves an integer number of terms,
so this can be evaluated classically. All the problems sit in the
last sum on the right-hand side. The payoff is that we have
translated the domain of summation by $n$ to the right. Since (\ref{EqHeuristics})
holds for every integer $n$, we can use $\axright$ to evaluate the
limit as $n\to\infty$: if $f(n+z)\to 0$ as $n\to\infty$ for all $z$,
then $\axright$ implies that the limit
as $n\to\infty$ of the last sum should vanish. We get
\[
   \sum_{\nu=x}^y f(\nu)=
   \sum_{\nu=1}^\infty \left(\strut f(\nu+x-1)-f(\nu+y)\right).
\]

This is of course a special condition to impose on $f$, but the
same idea can be generalized. For example, if $f(\nu)=\ln \nu$,
then for $\nu\in [x,y]\subset\R^+$, the values $f(\nu+n)$ are
approximated well by the constant function $f(n)$, with an
error that tends to $0$ as $n\to\infty$: we say that $f=\ln$ is
``approximately constant.'' Using (S3),
\[
   \sum_{\nu=x}^y \ln(\nu+n)=
   \sum_{\nu=x}^y \ln n +
   \sum_{\nu=x}^y \left(\strut
      \ln(\nu+n)-\ln n
   \right)
\]
for every $n\in\N$. But by $\axright$, the last sum vanishes as
$n\to\infty$, while the first sum on the right-hand side has a
constant summand and is evaluated using Proposition~\ref{Prop1}.
Taking the limit $n\to\infty$ in (\ref{EqHeuristics}),
it follows by necessity that
\[
   \sum_{\nu=x}^y \ln \nu=
   \lim_{n\to\infty}\left( \sum_{\nu=1}^n
   \left(\ln(\nu+x-1)-\ln(\nu+y)\strut\right)
   +(y-x+1)\ln n \right).
\]

Before generalizing our definition further, we take courage by
observing that this interpolates the factorial function in the
classical way: we define
\[
   \prod_{\nu=x}^y f(\nu):=
   \exp\left(\sum_{\nu=x}^y \ln f(\nu)\right)
\]
and thus get
\begin{eqnarray}
   z!=\prod_{\nu=1}^z \nu&=&
   \lim_{n\to\infty}
   \exp\left(
   \sum_{\nu=1}^n \ln\left(\frac{\nu}{\nu+z}\right)
   +z\ln n
   \right)\nonumber\\
   &=&
   \lim_{n\to\infty}\left(
   n^z\cdot
   \prod_{\nu=1}^{n}
   \frac{\nu}{\nu+z}\strut\right)
   =\Gamma(z+1),
   \label{EqGamma}
\end{eqnarray}
using a well-known product representation of the $\Gamma$
function~\cite[6.1.2]{AS}.

It is now straightforward to use the heuristic calculation in (\ref{EqHeuristics})
together with Proposition~\ref{Prop1} and axiom $\axright$
to derive a general definition: all we need is that the value of $f(n+z)$ can
be approximated by some sequence of polynomials $p_n(n+z)$ of fixed degree for $n\to\infty$.

Some care is needed with the domains of definition: the example of the
logarithm shows that it is inconvenient to
restrict to functions which are defined on all of $\C$. All we
need is a domain of definition $U$ with the property that $z\in
U$ implies $z+1\in U$. This leads to the following (using the convention that the zero polynomial is the unique polynomial of degree $-\infty$).

\begin{definition}[Fractional Summable Functions]
\label{DefFracSummable} \lineclear
Let $U\subset \C$ and $\sigma\in\N\cup\{-\infty\}$.
A function $f\colon U\to\C$ will be called {\em fractional summable
of degree $\sigma$} if the following conditions are satisfied:
\begin{itemize}
\item
$x+1\in U$ for all $x\in U$;
\item
there exists a sequence of polynomials $(p_n)_{n\in\N}$ of fixed
degree $\sigma$ such that for all $x\in U$
\[
\left| f(n+x) - p_n(n+x) \right| \longrightarrow 0
\quad \mbox{as $n\to+\infty$} \,;
\]
\item
for every $x,y+1\in U$, the limit
\[
    \lim_{n\to\infty}\left(
       \sum_{\nu=n+x}^{n+y}p_n(\nu)+\sum_{\nu=1}^n
   \left(\strut
          f(\nu+x-1)-f(\nu+y)
       \right)
    \right)
\]
exists, where $\sum p_n$ is defined as in~(\ref{eqSumPoly}).
\end{itemize}
In this case, we will use the notation
\[
\sumr_{\nu=x}^y f(\nu)
\qquad\mbox{or briefly}\qquad
\sumr_x^y f
\]
for this limit. Moreover, we can define fractional products by
\[
    \prodr_{\nu=x}^y f(\nu):=\exp\left(
       \sumr_{\nu=x}^y \ln f(\nu)
    \right),
\]
whenever $\ln f$ is fractional summable.
\end{definition}

Note that this definition does not depend on the choice of the approximating
polynomials $(p_n)_{n\in\N}$: if $(\tilde p_n)_{n\in\N}$ is another choice
of approximating polynomials, then $\lim_{n\to\infty} \left( p_n(n+x)-\tilde p_n(n+x)\right)=0$
for all $x\in U$, and hence for all $x\in\C$ since the set of polynomials of degree at most $\sigma$
is a finite-dimensional linear space.
As shown in Proposition~\ref{Prop1}, sums of polynomials
satisfy axiom $\axright$. Substituting $0$ for $f$ and $\tilde p_n-p_n$ for $p_n$ in $\axright$
proves that $\lim_{n\to\infty}\left(\sum_{\nu=n+x}^{n+y} p_n(\nu)-\sum_{\nu=n+x}^{n+y} \tilde p_n(\nu)\right)=0$.

Moreover, this definition is the unique definition that satisfies axioms (S1) to $\axright$:
\begin{theorem}
Definition~\ref{DefFracSummable} satisfies all the axioms (S1) to $\axright$ (for suitable domains of definition), and it is the unique definition with this property (for the class of functions that we are considering).
\end{theorem}
\begin{proof}
We have already proved uniqueness above, by deriving Definition~\ref{DefFracSummable}
from the axioms (S1) to $\axright$. It remains to prove that this definition indeed satisfies
all the axioms. Clearly, (S3) and (S5) are automatically satisfied. Substituting the
definition into (S1), (S2), and (S4), these axioms can be confirmed by a few lines of
direct calculation. To prove $\axright$, we use the definition and the other axioms
(in particular continued summation (S1)) and calculate
\begin{eqnarray*}
   \Delta&:=&\lim_{n\to\infty}\left(\sumr_{\nu=x}^y f(\nu+n)-\sum_{\nu=x}^y p_n(\nu+n)\right)\\
   &=&\lim_{n\to\infty}\left(
      \sumr_{\nu=x}^y f(\nu+n) - \sumr_{\nu=x}^y f(\nu) + \sum_{\nu=1}^n \left(f(\nu+x-1)-f(\nu+y)\right)
   \right)\\
   &=&\lim_{n\to\infty}\left(
      \sumr_{\nu=x+n}^{y+n} f(\nu)-\sumr_{\nu=x}^y f(\nu) + \sum_{\nu=x}^{n+x-1} f(\nu) - \sum_{\nu=y+1}^{y+n} f(\nu)
   \right)=0.
\end{eqnarray*}
This proves that Definition~\ref{DefFracSummable} satisfies all the axioms.
\end{proof}

\section{Properties of Fractional Sums.}
\label{SecProperties}
Now that we have a definition of sums with noninteger numbers of
terms, it is interesting to find out how many of the properties
of classical finite sums remain valid in this more general setting,
and what new properties arise that are not visible in the classical case.

\subsection{Generalized Classical Properties.}
One of the most basic identities for finite sums is
the geometric series. For simplicity, let $0\leq q <1$. Then
the function $\nu\mapsto q^\nu$ is approximately zero (we have $\lim_{n\to\infty}
q^{z+n}=0$ for every $z\in\C$), and the definition reads
\begin{equation}
   \sumr_{\nu=0}^x q^\nu = \sum_{\nu=1}^\infty \left(
      q^{\nu-1}-q^{\nu+x}
   \right)=\left( 1-q^{x+1}\right)\sum_{\nu=1}^\infty q^{\nu-1}=\frac{1-q^{x+1}}{1-q}.
   \label{EqGeometric}
\end{equation}
Thus, the formula for the geometric series remains valid for every $x\in\C$.

A similar calculation shows that the binomial series remains valid in the fractional case:
for every $c\in\C\setminus\{-1,-2,-3,\ldots\}$ and $x\in\C$ with $|x|<1$, we have
\begin{equation}
   (1+x)^c=\sumr_{\nu=0}^c {c \choose \nu} x^\nu.
   \label{EqBinomial}
\end{equation}

There are generalizations of (\ref{EqGeometric}) to the case $q>1$ and of (\ref{EqBinomial}) to the case $|x|>1$: these involve a ``left sum'' as introduced in Section~\ref{SecLeftSummation}.

An example of a summation identity with more complicated structure is given
by the series multiplication formula
\begin{eqnarray}
   \qquad\left(
      \sumr_{\nu=1}^x f(\nu)
   \right)&\cdot& \left(
      \sumr_{\nu=1}^x g(\nu)
   \right)\nonumber\\
   &=&\sumr_{\nu=1}^x\left(
      f(\nu)g(\nu)+f(\nu)\sumr_{k=1}^{\nu-1}g(k)
      +g(\nu)\sumr_{k=1}^{\nu-1} f(k)
   \right)
   \label{EqSeriesMultiplication}
\end{eqnarray}
for every $x\in\C$, given that all the three fractional sums exist (see \cite[Lemma~7]{FracSums}; it generalizes the formula $(a_1+a_2)(b_1+b_2)=a_1b_1+a_2b_2+a_2b_1+b_2a_1$, and similarly for all positive integers $x$.

\subsection{New Properties and Special Functions.}
As shown in Section~\ref{SecFromAxToDef}, our definition
interpolates the factorial by the $\Gamma$ function,
\begin{equation}
   z!\equiv \prodr_{\nu=1}^z \nu =\Gamma(z+1).
\end{equation}

An amusing consequence is
\begin{equation}
\prodr_{\nu=1}^{-1/2} (\nu^2+1) = \tanh \pi ;
\end{equation}
this is because 
\begin{eqnarray*}
\prodr_{\nu=1}^{-1/2}(\nu^2+1)&=&\prodr_{\nu=1}^{-1/2}(\nu+i)\prodr_{\nu=1}^{-1/2}(\nu-i)=
\frac{\Gamma(1/2+i)\Gamma(1/2-i)}{\Gamma(1+i)\Gamma(1-i)}\\
&=&\frac{\Gamma(1/2+i)\Gamma(1/2-i)}{i\Gamma(i)\Gamma(1-i)}
=\frac{\sin(\pi i)}{i\sin(\pi(1/2+i))} 
=\frac{\sin\,i\pi}{i\cos\,i\pi}\\
&=&\frac{\sinh\pi}{\cosh\pi}=\tanh\pi
\end{eqnarray*}
using $\Gamma(z)\Gamma(1-z)=\pi/\sin(\pi z)=i\pi/\sinh(\pi iz)$.

Many basic fractional sums are related to special
functions. As a first example, consider the harmonic series. Since $\nu\mapsto \nu^{-1}$
is approximately zero, the definition reads
\begin{equation}
   \sumr_{\nu=1}^x \frac 1 \nu=\sum_{\nu=1}^\infty\left(
      \frac 1 \nu - \frac 1 {\nu+x}
   \right),
\end{equation}
and in particular 
\begin{equation}
   \sumr_{\nu=1}^{-1/2}\frac 1 \nu= -2\left(1-\frac 1 2 +\frac 1 3-\frac 1 4+\cdots\right)=  -2\ln 2,
\label{Eq:Euler}
\end{equation}
which was noticed already by Euler~\cite[pp.\ 88--119, \S 19]{Euler}.
For general $x$, the harmonic series can be expressed in terms of the so-called digamma function
~\cite[6.3.1]{AS} $\psi(x+1)=\frac d {dx} \ln\Gamma(x+1)$ and the Euler-Mascheroni constant $\gamma=0.577\ldots$: one
obtains~\cite[6.3.16]{AS}
\begin{equation}
   \sumr_{\nu=1}^x \frac 1 \nu=\sum_{\nu=1}^\infty\left(
      \frac 1 \nu - \frac 1 {\nu+x}
   \right)
=\gamma+\psi(x+1).
\end{equation}

Note that the reflection formula~\cite[6.3.7]{AS} for the digamma function
becomes
\begin{equation}
   \sumr_{\nu=x}^{-x} \frac 1 \nu =\pi\cot(\pi x).
   \label{EqReflection}
\end{equation}

As a further generalization, it is convenient to consider the Hurwitz $\zeta$ function,
traditionally defined by the series
\[
   \zeta(s,x):=\sum_{\nu=0}^\infty \frac 1 {(\nu+x)^s},\qquad \Re(s)>1.
\]
By analytic continuation, $\zeta(s,x)$ can be defined for every $s\in\C$, except
for a pole at $s=1$. For $x=1$, the Hurwitz $\zeta$ function equals the well-known
Riemann $\zeta$ function:
\[
   \zeta(s,1)=\zeta(s).
\]
It turns out that the Hurwitz $\zeta$ function can be understood as
a fractional power sum. It can be shown \cite[Corollary~14]{FracSums} that
for every $a\in\C\setminus\{-1\}$ and for all $x\in\C\setminus\{-1,-2,-3,\ldots\}$,
\begin{equation}
   \sumr_{\nu=1}^x \nu^a=\zeta(-a)-\zeta(-a,x+1).
   \label{EqHurwitz}
\end{equation}

A useful special case is
\begin{equation}
   \sumr_{\nu=1}^{-\frac 1 2} \nu^a =\left(2-2^{-a}\right)\zeta(-a).
   \label{EqZeta12}
\end{equation}
Note that such equations give in many cases intuitive ways to compute
properties and special values of special functions. Everybody knows the
formula $\sum_{\nu=1}^x \nu =x(x+1)/2$, so $\sumrinline_{\nu=1}^{-1/2}\nu=-1/8$, and thus by (\ref{EqZeta12})
\[
   -\frac 1 8 =\sumr_{\nu=1}^{-1/2} \nu^1=\left(2-2^{-1}\right)\zeta(-1)=\frac 3 2 \zeta(-1).
\]
It follows that $\zeta(-1)=-1/12$. Similarly, we have
\[
   \frac d {dz}\left(2-2^{-z}\right) \zeta(-z)\Big\vert_{z=0}
   =\frac d {dz} \left.\sumr_{\nu=1}^{-1/2} \nu^z \right|_{z=0}=\sumr_{\nu=1}^{-1/2} \ln\nu
   =\ln\prodr_{\nu=1}^{-1/2}\nu=\ln\Gamma\left(\frac 1 2\right)
\]
(in the second equality, we interchanged differentiation and fractional summation; it is not hard to check that this is indeed allowed).
Since $\Gamma(1/2)=\sqrt{\pi}$, this easily implies that $\zeta'(0)=-(1/2) \ln(2\pi)$.

Similarly, differentiating (\ref{EqHurwitz}) $b$ times with respect to $a$ and arguing as before (compare also \cite[Sec.~6]{FracSums}), we obtain 
\begin{equation}
   \sumr_{\nu=1}^x \nu^a (\ln\nu)^b=(-1)^b\left(
      \zeta^{(b)}(-a)-\zeta^{(b)}(-a,x+1)
   \right).
   \label{EqHurwitzDiff}
\end{equation}

There are some classically unexpected special values like
\begin{equation}
   \sumr_{\nu=1}^{-1/2}\nu\ln\nu=-\frac{\ln 2}{24} -\frac 3 2 \zeta'(-1).
   \label{Eqvlnv}
\end{equation}

\subsection{Mirror Series and Left Summation.}
There is an identity for classical sums which is almost never mentioned,
because it seems so trivial. Consider the sum
\[
   f(-10)+f(-9)+f(-8)+f(-7).
\]
Obviously, there are two formally correct possibilities to write
this sum,
\[
   \mbox{either }\sum_{\nu=-10}	^{-7} f(\nu)\quad\mbox{ or }\quad
   \sum_{\nu=7}^{10} f(-\nu).
\]
Classically, it is clear that $\sum_{\nu=a}^b f(\nu)=\sum_{\nu=-b}^{-a} f(-\nu)$.
Does this carry over to the fractional case? There is a fundamental problem: our definition of fractional sums involves $\lim_{n\to+\infty}f(\nu+n)$, i.e., $f$ is evaluated near $+\infty$, and when $f(\nu)$ is replaced by $f(-\nu)$ then $f$ would be evaluated near $-\infty$ where the values may be unrelated. This will be discussed in the next section.

\section{An Alternative Axiom and Left Summation.}
\label{SecLeftSummation}
Looking back at the axioms given in Section~\ref{SecTheAxioms}, there is one
axiom that could possibly be modified: in $\axright$,
limits as $n\to +\infty$ are considered, but one could equally well look at
limits as $n\to -\infty$. This way, one obtains an axiom of ``left shift
continuity'':

\begin{description}
\item[\boldmath$\axleft$ \unboldmath Left Shift Continuity]
if $\lim_{n\to\infty}f(z-n)=0$ pointwise for $z\in\C$, then
\begin{equation}
   \lim_{n\to\infty} \sum_{\nu=x}^y f(\nu-n)=0;
   \label{EqS5L}
\end{equation}
   more generally, if there is a sequence of polynomials $(p_n)_{n\in\N}$
   of fixed degree such that $|f(z-n)-p_n(z-n)|\longrightarrow 0$ for $z\in\C$ as $n\to\infty$,
   we require that
   \[
      \left\vert \sum_{\nu=x}^y f(\nu-n)-\sum_{\nu=x}^y p_n(\nu-n)\right\vert\longrightarrow 0.
   \]
\end{description}

Repeating the calculations of Section~\ref{SecFromAxToDef}, one gets an alternative
definition\footnote{Note that no other complex directed
limit to infinity (like $n\to i\infty$) can determine a definition uniquely: only adding or
subtracting $n\in\N$ to the upper summation boundary consists of adding or subtracting
$n$ terms to the series, which can be done classically.}
which we do not state here formally: it is exactly
the same as Definition~\ref{DefFracSummable}, except that in every
limit, $n\to\infty$ is replaced by $n\to -\infty$.

It can be shown that this definition is the unique one that satisfies axioms (S1), (S2), (S3),
(S4), (S5), and $\axleft$. Note that in general,
the existence of $\sumr$ and $\suml$ are independent, and if both left and right fractional sums exist, they may
have different values. For example, for every $z$ with $\Re(z)>0$, we have
\[
   \suml_{\nu=1}^{-1/2}\nu^z=(-1)^{z+1} \left(2-2^{-z}\right)\zeta(-z),
\]
in contrast to equation~(\ref{EqZeta12}).

What we do have is the obvious relation
   \begin{equation}
      \sumr_{\nu=a}^b f(\nu)=\suml_{\nu=-b}^{-a}f(-\nu).
\label{Eq:MirrorSeries}   
\end{equation}

\section{Classical Infinite Sums, Products, and Limits.}
\label{SecClassical}

Fractional sums are not simply a new world with results that have
no meaning in the classical context; they allow us to derive identities
that can be stated entirely in classical terms. Some of these
formulas are known and some seem to be new. Of course, all these
identities can in principle be computed without fractional sums.
But proving them with the help of fractional sums is rather
intuitive and simple, since most of the steps use fractional generalizations
of basic, very well-known classical summation properties.

\subsection{Some Infinite Products.}
As a first example, we show how to compute a closed-form expression
for the infinite product
\begin{equation}
   P(x):=\lim_{n\to\infty} \prod_{k=1}^{2n} \left(1+\frac {2x} k\right)^{-k(-1)^k}
   \label{eq:defPx}
\end{equation}
for $x>-1/2$. It was first considered by Borwein and Dykshoorn in 1993 (see \cite{BD}). By
taking logarithms, one gets
\[
   \ln P(x)=-\sum_{k=1}^\infty \left(
      2k\ln\left(1+\frac {2x}{2k}\right)-2\left(k-\frac 1 2\right)\ln\left(1+\frac{2x}{2\left(k-1/2\right)}\right)
   \right).
\]
Consider the function $\nu\mapsto 2\nu\ln\left(1+x/\nu\right)$ that tends
to $2x$ as $\nu\to\infty$. According to Definition~\ref{DefFracSummable}, we have
\begin{eqnarray*}
   &&\sumr_{\nu=1}^{-1/2}2\nu\ln\left(1+\frac x \nu\right)\\
   &=&\lim_{n\to\infty}\left[ 
       -\frac 1 2 \cdot 2x 
        +\sum_{k=1}^n 2 k \ln\left(1+\frac x k\right)
      -2\left(k-\frac 1 2\right)\ln\left(1+\frac x {k-1/2}\right)
   \right].
\end{eqnarray*}
Thus, we get
\begin{eqnarray*}
   \ln P(x)&=& -x-\sumr_{\nu=1}^{-1/2} 2\nu\ln\left(1+\frac x \nu\right)=
   -x-2\sumr_{\nu=1}^{-1/2} \nu\ln\left(\frac{\nu+x}{\nu}\right)\\
   &=&-x-2\sumr_{\nu=1}^{-1/2}\nu\ln(\nu+x)+2\sumr_{\nu=1}^{-1/2}\nu\ln\nu\\
   &=&-x-2\sumr_{\nu=1+x}^{-1/2 +x} (\nu-x)\ln\nu -\frac {\ln 2}{12} -3 \zeta'(-1)\\
   &=&-x-2\sumr_{\nu=1+x}^{-1/2 +x} \nu\ln\nu +2 x \sumr_{\nu=1+x}^{-1/2 +x} \ln\nu
   -\frac {\ln 2}{12} -3 \zeta'(-1)\\
   &=&-x-2\left(\zeta'\left(-1,x+\frac 1 2\right)-\zeta'\left(-1,x+1\right)\right)\\
   && + 2 x \left(\ln\left(\left(x-\frac 1 2 \right)!\right)-\ln(x!)\right)-\frac {\ln 2}{12} -3 \zeta'(-1),
\end{eqnarray*}
where we have used equation~(\ref{Eqvlnv}), index shifting, equation~(\ref{EqHurwitzDiff}),
equation~(\ref{EqGamma}), and continued summation. By exponentiating, we finally get
\begin{equation}
   P(x)=2^{-\frac 1 {12}} \left(
      \frac{\Gamma\left(x+\frac 1 2\right)}{\Gamma(x+1)}
   \right)^{2 x} e^{
      -x-2\zeta'\left(-1,x+\frac 1 2\right)+2\zeta'\left(-1,x+1\right)-3\zeta'(-1)
   }.
   \label{eq:BD}
\end{equation}

\begin{figure}
\centerline{
\includegraphics[width=0.6\textwidth]{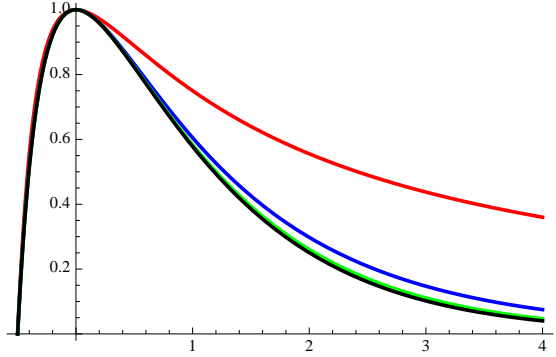}
}
\caption{Numerical check of (\ref{eq:BD}). The right-hand side corresponds to the lowest curve,
while the other three curves (from top to bottom) are plots of approximations to $P(x)$ (finite products
as in~(\ref{eq:defPx})) for $n=1$, $n=10$, and
$n=50$ respectively.
}
\label{fig:BD}
\end{figure}

Using Mathematica's built-in numerical procedures, this infinite product identity can be checked numerically.
Figure~\ref{fig:BD} shows a comparison of both sides of this equation.

By application of this method, a large class of infinite products can be
explicitly computed, which seems to include the class of products
considered in \cite{Adamchik}. Here is an example of a new identity: using
the same steps as in the calculation above, one easily proves that for $x>-1$,
\begin{eqnarray*}
   \sumr_{k=1}^{-1/2} 2 k \ln^2(2k+x)
   &=&4\ln 2\left(\zeta'\left(-1,\frac{x+1} 2\right)
   -\zeta'\left(1,\frac x 2 +1\right)\right)\\
   &&-2 x \ln 2 \ln\frac{\Gamma\left(\frac{x+1} 2\right)}{\Gamma\left(\frac x 2 +1\right)}
   +2\zeta''\left(-1,\frac x 2+1\right)\\
   &&-2\zeta''\left(-1,\frac{x+1} 2\right)+x\zeta''\left(0,\frac{x+1} 2\right)\\
   &&-x\zeta''\left(0,\frac x 2+1\right)-\frac {\ln^2 2} 4.
\end{eqnarray*}

\begin{figure}
\centerline{
\includegraphics[width=0.6\textwidth]{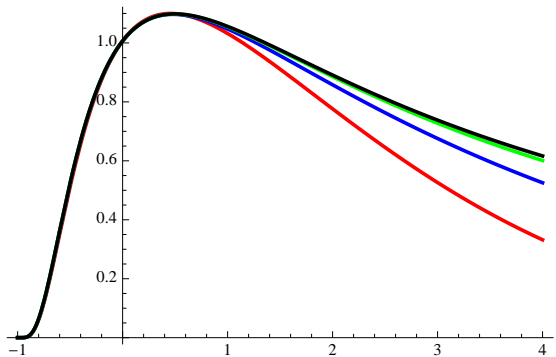}
}
\caption{Numerical check of (\ref{EqComplicated}). The right-hand side corresponds to the uppermost curve,
while the other three curves (from bottom to top) are plots of the left-hand side for $n=10$, $n=100$, and
$n=1000$ respectively.
}
\label{fig:secondzeta}
\end{figure}

Resolving the definition and exponentiating, we get the following
classical limit identity:
\begin{eqnarray}
&&\lim_{n\to\infty}\left[
(2n)^{-\frac 1 2 -x -\left(n+\frac 1 4\right)\ln(2n)}\prod_{k=1}^{2n}
(k+x)^{(-1)^k k \ln(k+x)}\right]\nonumber\\
&=& 2^{-\frac 1 4 \ln 2 +4 \zeta'\left(-1,\frac{x+1} 2\right)-4\zeta'\left(-1,\frac x 2 +1\right)}
\left(
   \frac{\Gamma\left(\frac{x+1} 2 \right)}
   {\Gamma\left(\frac x 2+1\right)}
\right)^{-2x\ln 2}\times\nonumber\\
&&\times\quad e^{
   2 \zeta''\left(-1,\frac x 2+1\right)-2\zeta''\left(-1,\frac{x+1} 2\right)
   +x\left(
      \zeta''\left(0,\frac {x+1} 2\right)-\zeta''\left(0,\frac x 2+1\right)
   \right)
}.\label{EqComplicated}
\end{eqnarray}
Again, we have used Mathematica for a quick numerical check that is shown in Figure~\ref{fig:secondzeta}.

\subsection{The Multiple $\Gamma$ Function and Series Multiplication.}
In this section, we consider the multiple gamma function $\Gamma_n$, a
generalization of the classical gamma function $\Gamma$, defined
for $n\in\N$ and $z\in\C$ by the recurrence formula (compare \cite{Adamchik})
\begin{eqnarray}
   \Gamma_{n+1}(z+1)&=&\frac {\Gamma_{n+1}(z)}{\Gamma_n(z)},\nonumber\\
   \Gamma_0(z)&=&z^{-1},\label{EqDefGammaN}\\
   \Gamma_n(1)&=&1.\nonumber
\end{eqnarray}
These equations do not determine the functions $\Gamma_n$ uniquely, so one
needs the additional Bohr-Mollerup-like condition that $\Gamma_n(x)$ is positive
and $n$ times differentiable on $x\in\R^+$, and that $(-1)^{n+1}\frac{d^n}{dx^n} \ln \Gamma_n(x)$ is increasing
(see \cite{Duke}). For $n=1$, this definition reproduces the classical gamma function: $\Gamma_1(z)=\Gamma(z)$.

The (reciprocal of the) special case $n=2$ is known as the Barnes $G$ function
\[
   G(z)=1/\Gamma_2(z).
\]
By (\ref{EqDefGammaN}), it satisfies
\[
   G(n)=\Gamma(1)\Gamma(2)\Gamma(3)\cdots\Gamma(n-1).
\]
More generally,
\[
   1/\Gamma_{n+1}(z)=\Gamma_n(1)\Gamma_n(2)\Gamma_n(3)\cdots\Gamma_n(z-1),
\]
so $\Gamma_{n+1}$ is the reciprocal of the product of $\Gamma_n$, which means
that $\Gamma_n(z)$ is something like an $n$-fold product of the first $z-1$ natural numbers:
\[
   \left(\Gamma_n(z)\right)^{(-1)^{n+1}}=\prod_{\nu_1=1}^{z-1}\prod_{\nu_2=1}^{\nu_1-1}\prod_{\nu_3=1}^{\nu_2-1}
   \cdots\prod_{\nu_n=1}^{\nu_{n-1}-1}\nu_n.
\]
While this equation only makes sense for $z\in\N\setminus\{0\}$, one can easily show that the definition
of $\Gamma_n(z)$ for $z\in\C$ is compatible with our definition for fractional sums and
products, i.e., that for every $z\in\C$ (except for poles) and $n\in\N\setminus\{0\}$, we have
\[
   (-1)^{n+1} \ln \Gamma_n(z)=\sumr_{\nu_1=1}^{z-1}\sumr_{\nu_2=1}^{\nu_1 -1}\sumr_{\nu_3=1}^{\nu_2-1}
   \ldots\sumr_{\nu_n=1}^{\nu_{n-1}-1} \ln \nu_n.
\]
We will now show some properties of the multiple gamma function $\Gamma_n$,
specifically for the example $n=2$, simply by using basic fractional sum identities,
without using any special function properties of $\Gamma_n$.
By the multiplication formula (\ref{EqSeriesMultiplication}), we have
\begin{eqnarray*}
   \ln G(z)&=&-\ln \Gamma_2(z)=\sumr_{\nu=1}^{z-1}\left(1\cdot \sumr_{k=1}^{\nu-1} \ln k\right)\\
   &=&\left(\sumr_{\nu=1}^{z-1} 1\right)\left(\sumr_{\nu=1}^{z-1}\ln\nu\right)
   -\sumr_{\nu=1}^{z-1}1\cdot\ln\nu-\sumr_{\nu=1}^{z-1}\left(\ln\nu\sumr_{k=1}^{\nu-1} 1\right)\\
   &=&(z-1)\ln\Gamma(z)-\sumr_{\nu=1}^{z-1}\nu\ln\nu\\
   &=&(z-1)\ln\Gamma(z)+\zeta'(-1)-\zeta'(-1,z).
\end{eqnarray*}
The last equality follows from equation~(\ref{EqHurwitzDiff}). Thus, we have found
an explicit formula for $G(z)$ in terms of derivatives of the Hurwitz $\zeta$ function.

Equation~(\ref{Eqvlnv}) gives the special value
\[
   G\left(\frac 1 2\right)=e^{
      -\frac 1 2 \ln\Gamma\left(\frac 1 2\right)-\sumrinline_{\nu=1}^{-1/2} \nu\ln\nu
   }
   =\pi^{-\frac 1 4}\,\,2^{\frac 1 {24}}\,\,e^{\frac 3 2 \zeta'(-1)}.
\]
These are of course very well-known results, but the calculations are strikingly simple.
Moreover, this example shows that there is a wide variety of interesting ``special functions''
that do not have to be defined separately, but can be treated in a unified manner
by our theory of fractional sums.
New generalizations comparable to $G(z)=\prodr_{n=0}^{z-2} n!$ include
\begin{eqnarray*}
   \prodr_{n=1}^{-1/2} (2n)!&=&\left(\frac \pi 2\right)^{\frac 1 4},\\
   \prodr_{n=1}^{-1/2}\left(n!\right)^{\ln n}&=&\exp\left(\frac{\gamma^2}4+\frac{\gamma_1} 2
   -\frac{\pi^2}{48}+\frac{\ln^2 2}2-\frac{\ln^2\pi}8\right),
\end{eqnarray*}
\[
   \prodr_{n=1/4}^{-1/4} \left(n!\right)^n
   =\left(\frac{\Gamma\left(\frac 1 4\right)}{\Gamma\left(\frac 34\right)}\right)^{\frac 3 {32}}
   e^{
      \zeta'\left(-2,\frac 1 4\right)-\frac{3\zeta(3)}{128\pi^2} -\frac G {4\pi}
   }.
\]
Here, $\gamma=0.577215\ldots$, $\gamma_1=.072815\ldots$, and $G=.91596\ldots$ are the Euler-Mascheroni,
Stieltjes, and Catalan constants, respectively.
Again, these formulas have classical limit representations looking like equation~(\ref{EqComplicated})
which we do not write down here explicitly.

\subsection{Perspective: A Series by Gosper.}
The paper ``On some strange summation formulas'' \cite{Strange} contains
some formulas like (\ref{EqGosper}) below. There might possibly be very short proofs
for all these identities using fractional sums. The only problem is that
there is one single step (indicated by the question mark) which we are unable
to justify: it is basically an interchange of a fractional
sum and an infinite series.

Nevertheless, we give this calculation as a speculation, just to show that it
is tempting to have a closer look at what else might still be possible.

{\sc Speculation} (A Series by Gosper).
\lineclear
{\em For every $b\in\R$, we have the identity}
\begin{equation}
   S(b):=\sum_{n=0}^\infty \frac{(-1)^n}{(n+\frac 1 2)}\frac
   {\sin \sqrt{b^2+\pi^2 (n+1/2)^2}}
   {\sqrt{b^2+\pi^2 (n+1/2)^2}}
   =\frac{\pi \sin b}{2 b}.
   \label{EqGosper}
\end{equation}

{\sc ``Proof.''} We start by writing the aforementioned series as a
fractional sum with $-1/2$ terms. By plugging in the definition,
one easily confirms that
\begin{equation}
   \sum_{n=0}^\infty \frac{(-1)^n}{(n+\frac 1 2)}\frac
   {\sin \sqrt{b^2+\pi^2 (n+1/2)^2}}
   {\sqrt{b^2+\pi^2 (n+1/2)^2}}
   =-\sumr_{n=3/4}^{-3/4} \frac 1 {2n}
   \frac{\sin\sqrt{b^2+4 \pi^2 n^2}}{\sqrt {b^2+4\pi^2 n^2}}.
   \label{EqSb}
\end{equation}

We will now use the basic identity
\begin{equation}
   \sumr_{\nu=x}^{-x} \nu^{2 n+1}=0
   \label{EqSumxminusx}
\end{equation}
for every $x\in\C$ and $n\in\N$, which can be shown in two different ways:
The first possibility is to see that for every $x\in -\N$,
\[
   \sumr_{\nu=x}^{-x} \nu^{2 n+1}=x^{2n+1}+(x+1)^{2n+1}+\cdots+(-x-1)^{2n+1}+(-x)^{2n+1}=0,
\]
since even and odd terms cancel each other. By continued summation and Proposition~\ref{Prop1},
$\sumrinline_{\nu=x}^{-x} \nu^{2 n+1}=
\sumrinline_{\nu=1}^{-x}\nu^{2 n+1} - \sumrinline_{\nu=1}^{x-1}\nu^{2 n+1}$ is a polynomial in $x$,
so equation~(\ref{EqSumxminusx}) must be valid for every $x\in\C$.

A second way is to use the mirror series from (\ref{Eq:MirrorSeries}) to calculate
\begin{equation}
   \sumr_{\nu=x}^{-x} \nu^{2 n+1}=\suml_{\nu=x}^{-x} (-\nu)^{2 n+1}=-\suml_{\nu=x}^{-x} \nu^{2 n+1}.
   \label{EqMirrorInAction}
\end{equation}
For polynomials, left and right sum coincide trivially, so (\ref{EqSumxminusx}) follows immediately.

Going back to the fractional sum in (\ref{EqSb}), the odd function
\begin{equation}
   f(n):=\frac 1 {2n}
   \frac{\sin\sqrt{b^2+4 \pi^2 n^2}}{\sqrt {b^2+4\pi^2 n^2}}
   \label{Eqfn}
\end{equation}
is holomorphic in the entire complex plane, except for a pole at $n=0$, so we can develop it into a power series. We get
\[
   S(b)=- \sumr_{n=3/4}^{-3/4}\left(
      \frac{\sin b}{2 b}n^{-1}+c_1 n +c_3 n^3 +c_5 n^5 +\cdots
   \right).
\]
The next step is critical: we apply the fractional sum term-by-term.
Unfortunately, it is not clear that this manipulation is justified.
\begin{equation}
   S(b)\stackrel{?}{=}-\left(
      \frac{\sin b}{2 b}\sumr_{n=3/4}^{-3/4} n^{-1}
      +c_1 \sumr_{n=3/4}^{-3/4} n +c_3 \sumr_{n=3/4}^{-3/4} n^3+\ldots
   \right).
   \label{EqTermwise}
\end{equation}
Equations~(\ref{EqReflection}) and (\ref{EqSumxminusx}) yield
\[
   S(b)=-\frac{\sin b}{2b}\pi\cot\left(\frac 3 4 \pi\right)=\frac{\pi \sin b}{2 b}.
\]
\qed

This method only works for a certain class of functions which obviously contains $f(n)$ from (\ref{Eqfn})
and other functions like $e^{an}/n^k$, but which does {\em not} contain other simple functions
like $e^{-a n^2}$. It is an open question to give sufficient conditions for the validity
of this method, i.e., for justification of termwise fractional summation as in
equation~(\ref{EqTermwise}).

\paragraph{Acknowledgments.}  We would like to thank several colleagues from the community of ``special functions and exotic identities'' for their encouragement and support, especially Richard Askey and Mourad Ismail. Moreover, we are grateful to Otto Forster, Irwin Kra, Armin Leutbecher,
John Milnor, as well as to the seminar audiences in Stony Brook and M\"unchen for encouragement, interest, and helpful discussions.
We would also like to thank the referees for helping us improve the exposition.

\bigskip

\noindent\textbf{Markus M\"uller} received his Dr.\ rer.\ nat.\ from Technical University of Berlin in 2007,
where he is currently working as a postdoc in quantum information theory.
After playing around with fractional sums as a high-school student, he was lured
away from math to physics by popular scientific articles on Schr\"odinger's cat and quantum weirdness.
Since then, his main motivation has been the idea that the notion of information opens up an
unexpected, fresh perspective on foundational problems of physics. This line of thought
led him to work on quantum Turing machines and Kolmogorov complexity, concentration of measure,
and generalized probabilistic theories beyond quantum theory.
At the moment, he is preparing for some tough postdoc years abroad
by skiing with friends and enjoying the alternative music occasions in Berlin.

\noindent\textit{Institute of Mathematics, Technical University of Berlin, Stra\ss e des 17.\ Juni 136, 10623 Berlin, Germany.}

\bigskip

\noindent\textbf{Dierk Schleicher} studied physics and computer science in Hamburg and obtained his Ph.D. in mathematics at Cornell University. He enjoyed longer educational and research visits in Princeton, Berkeley, Stony Brook, Paris, and Toronto. After many years in M\"unchen, he became the first professor at the newly-founded Jacobs University Bremen in 2001 and built up the mathematics program there. His main research area is dynamical systems, especially complex dynamics: ``real mathematics is difficult, complex mathematics is beautiful.'' He has always been active in math circles and special programs for talented high school students, which is what lured him away from physics to mathematics. He was one of the main organizers of the 50th International Mathematical Olympiad (IMO) in Bremen/Germany, in 2009. He enjoys outdoor activities: kayaking, paragliding, mountain hiking, and more.

\noindent\textit{Jacobs University Bremen, Research I, Postfach 750 561, D-28725 Bremen, Germany.}

\end{document}